\documentclass[graybox]{svmult}

\usepackage{type1cm}        % activate if the above 3 fonts are
\usepackage{makeidx}         % allows index generation
\usepackage{graphicx}        % standard LaTeX graphics tool
\usepackage{multicol}        % used for the two-column index
\usepackage[bottom]{footmisc}% places footnotes at page bottom

\usepackage{newtxtext}       % 
\usepackage[varvw]{newtxmath}       % selects Times Roman as basic font

\usepackage{cleveref}
\usepackage{subcaption}
\usepackage{easyReview}
\usepackage{ulem}

\DeclareMathOperator{\diverg}{div} % div operator

\makeindex             % used for the subject index
\begin{document}

\title*{Two-level additive Schwarz preconditioners for reduced integration methods}
\author{Filipe Cumaru\orcidID{0009-0003-9516-4226}, Alexander Heinlein\orcidID{0000-0003-1578-8104} and Joachim Sch\"oberl\orcidID{0000-0002-1250-5087}}
\institute{Filipe Cumaru \at Delft Institute of Applied Mathematics, Faculty of Electrical Engineering, Mathematics and Computer Science, Delft University of Technology, Delft, The Netherlands, \email{f.a.cumarusilvaalves@tudelft.nl}
\and Alexander Heinlein \at Delft Institute of Applied Mathematics, Faculty of Electrical Engineering, Mathematics and Computer Science, Delft University of Technology, Delft, The Netherlands, \email{a.heinlein@tudelft.nl}
\and Joachim Sch\"oberl \at Institute for Analysis and Scientific Computing, Vienna University of Technology, Vienna, Austria, \email{joachim.schoeberl@tuwien.ac.at}}
\maketitle

\abstract*{Incompressible fluid flow problems appear frequently in different applications. The discretization of such problems may result in large and ill-conditioned systems of linear equations. We consider the case of the Stokes equations discretized using a reduced integration method which approximates the incompressibility constraint by a penalty term thus allowing the problem to be solved only in terms of the velocity unknowns. We investigate the numerical scalability of a two-level overlapping additive Schwarz method with a reduced dimension generalized Dryja-Smith-Widlund (RGDSW) coarse space. In addition, we discuss the parallel implementation of the examples using the Fast and Robust Overlapping Schwarz (FROSch) package for additive Schwarz preconditioners and the NGSolve library, which implements multiple finite element space formulations.}

\abstract{Incompressible fluid flow problems appear frequently in different applications. The discretization of such problems may result in large and ill-conditioned systems of linear equations. We consider the case of the Stokes equations discretized using a reduced integration method which approximates the incompressibility constraint by a penalty term thus allowing the problem to be solved only in terms of the velocity unknowns. We investigate the numerical scalability of a two-level overlapping additive Schwarz method with a reduced dimension generalized Dryja-Smith-Widlund (RGDSW) coarse space. In addition, we discuss the parallel implementation of the examples using the Fast and Robust Overlapping Schwarz (FROSch) package for additive Schwarz preconditioners and the NGSolve library, which implements multiple finite element space formulations.}

% ------------

\section{Introduction} \label{sec:intro}

Incompressible fluid flow problems appear in different physical applications. The finite element (FE) discretization of such problems may result in large and ill-conditioned linear systems requiring effective preconditioning techniques to accelerate their iterative solution.

We consider a two-level overlapping additive Schwarz (OAS) preconditioner with a reduced dimension generalized Dryja-Smith-Widlund (RGDSW) coarse space \cite{DOHRMANN2017RGDSW}, a variant of the generalized Dryja-Smith-Widlund (GDSW) coarse space \cite{DOHRMANN2008GDSW1, DOHRMANN2008GDSW2} originally proposed for elliptic problems. The GDSW and the RGDSW coarse spaces have been successfully applied in a monolithic framework for a mixed formulation of incompressible flow problems \cite{HEINLEIN2019MONOLITHICGDSW,HEINLEIN2020MONOLITHICRGDSW}. Similar preconditioners with a Lagrangian coarse space \cite{KLAWONN2000COMPMONOOAS} as well as multigrid approaches \cite{SCHOEBERL1998,VANKA1986,VERFURTH1984} have been introduced.

In this work, we examine as a model problem the Stokes equations in three dimensions. We seek the velocity $u \in V_g = \{ v \in \left( H^1(\Omega) \right)^3 | v_{\partial \Omega_D} = g \}$ and pressure $p \in Q = L^2(\Omega)$ of an incompressible fluid such that
\begin{equation} \label{eq:weak_form}
    \begin{alignedat}{2}
        \int_{\Omega} \nabla u \colon \nabla v \; \text{d}x + \int_{\Omega} \diverg u \, p \; \text{d}x &= \int_{\Omega} f v \; \text{d}x && \quad \forall v \in V_0, \\
        \int_{\Omega} \diverg u \, q \; \text{d}x &= 0 && \quad \forall q \in L^2(\Omega).
    \end{alignedat}
\end{equation}

For the discretization of the problem in \cref{eq:weak_form} we employ a reduced integration method with a high order Lagrangian finite element space noted as $P_k$ for $k \geq 2$. The incompressibility constraint is approximated by a penalty term allowing us to solve the problem for the primal variable $u$ only thus reducing the size of the resulting linear system \cite{MALKUS1978,SCHOEBERL1998}. The modified bilinear form is given by
\begin{equation} \label{eq:bilinear_form_penalty}
    a_h(u, v) = \int_{\Omega} \nabla u \colon \nabla v \; \text{d}x + \epsilon^{-1} \int_{\Omega} \overline{\diverg u}^h \, \overline{\diverg v}^h \; \text{d}x,
\end{equation}
where $\epsilon \in (0, 1]$ is the penalty parameter and $\overline{\diverg u}^h$ represents an element-wise averaging of the divergence on the finite element mesh.

This article is structured as follows. In \cref{sec:two_level_precond}, we introduce the two-level OAS preconditioner with the RGDSW coarse space for the Stokes equations. Next, in \cref{sec:frosch_ngsolve}, we discuss the parallel implementation of the proposed preconditioner based on the FROSch package in Trilinos \cite{HEINLEIN2020FROSCH} and the NGSolve library \cite{NGSOLVE2014}. Finally, in \cref{sec:results}, we present numerical results to evaluate the parallel performance of the proposed preconditioning schemes.

% ------------

% ------------

\section{Two-level Schwarz preconditioners for saddle point problems} \label{sec:two_level_precond}

Let $A \mathbf{u} = b$ be the system of linear equations resulting from the discretization of \cref{eq:bilinear_form_penalty} on a computational domain $\Omega$. Furthermore, let $\Omega$ be decomposed into $N$ non-overlapping subdomains $\Omega_1,\ldots,\Omega_N$. Each non-overlapping subdomain is extended by layers of mesh elements to obtain an overlapping decomposition $\Omega_1',\ldots,\Omega_N'$. The corresponding two-level OAS preconditioner can be written as
\begin{equation} \label{eq:two_level_oas}
    M_{OAS, 2}^{-1} = \Phi A_0^{-1} \Phi^\top + \sum_{i = 1}^N R_i^\top A_i^{-1} R_i,
\end{equation}
where $R_i$ are the restriction matrices to the overlapping subdomains $\Omega_i'$, $A_i = R_i A R_i^\top$ are the corresponding local subdomain matrices, and $A_0 = \Phi^\top A \Phi$ is the coarse problem matrix. The coarse problem is defined on the coarse space spanned by a set of coarse basis functions which form the columns of the prolongation operator $\Phi$. In general, a coarse problem is required for the numerical scalability of domain decomposition methods; see, for instance~\cite{TOSELLI2005DDBOOK}.

\begin{figure}[t]
    \centering
    \begin{subfigure}[b]{0.3\textwidth}
        \includegraphics[width=\textwidth]{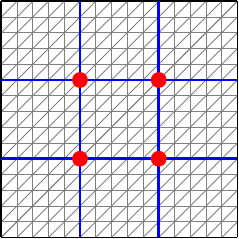}
        \caption{}
        \label{fig:dd_interface_2d}
    \end{subfigure}
    \hfil
    \begin{subfigure}[b]{0.3\textwidth}
        \includegraphics[width=\textwidth]{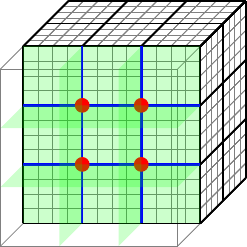}
        \caption{}
        \label{fig:dd_interface_3d}
    \end{subfigure}
    \caption{Interface of a non-overlapping domain decomposition in two (a) and three (b) dimensions. The vertex, edge and face interface components are marked in red, blue and green, respectively.}
    \label{fig:interface_components}
\end{figure}

In order to define the coarse space, let the interface of the non-overlapping domain decomposition be defined as
\begin{equation*} \label{eq:dd_interface}
    \Gamma = \bigcup_{i=1}^N \partial\Omega_i \setminus \partial\Omega_D.
\end{equation*}
We may thus reorder the degrees of freedom (dofs) in the system matrix $A$ and the prolongation operator $\Phi$ according to the division of the domain into interface ($\Gamma$) and interior (I $= \Omega \setminus \Gamma$):
\begin{equation} \label{eq:les_reorder}
    A = \begin{pmatrix}
        A_{\text{II}} & A_{\text{I}\Gamma} \\
        A_{\Gamma \text{I}} & A_{\Gamma \Gamma}
    \end{pmatrix}, \,
    \Phi = \begin{pmatrix}
        \Phi_{\text{I}} \\
        \Phi_{\Gamma}
    \end{pmatrix}.
\end{equation}
As usual, if the boundary dofs have not been eliminated, we treat them as interior. The value of the coarse basis functions in the interior $\Phi_{\text{I}}$ can be expressed as discrete harmonic extensions of the interface values $\Phi_{\Gamma}$
\begin{equation} \label{eq:discrete_harmonic_ext}
    \Phi = 
    \begin{pmatrix}
        -A_{\text{II}}^{-1} A_{\text{I} \Gamma} \\
        I_{\Gamma \Gamma}
    \end{pmatrix} \Phi_{\Gamma},
\end{equation}
where $I_{\Gamma \Gamma}$ is the identity matrix of dimensions $|\Gamma| \times |\Gamma|$. Multiple coarse spaces may be characterized based on different definitions of $\Phi_{\Gamma}$. This construction is employed in the definition of the GDSW coarse space \cite{DOHRMANN2008GDSW1,DOHRMANN2008GDSW2}.

\subsection{The RGDSW coarse space} \label{sec:rgdsw}

The reduced dimension generalized Dryja-Smith-Widlund (RGDSW) \cite{DOHRMANN2017RGDSW} coarse space is a variant of the GDSW \cite{DOHRMANN2008GDSW1, DOHRMANN2008GDSW2} coarse space with a smaller number of interface components used in the construction the coarse basis functions. While in GDSW, the interface $\Gamma$ is decomposed into vertex, edge and, in 3D, face components, in RGDSW, these components are grouped based on the subdomains that share them. An example of decomposition of the interface is illustrated in \cref{fig:interface_components}.

Let $\mathcal{N}_j$ be the set of nodes shared by the same subdomains, referred henceforth as a nodal equivalence class (nec). In addition, let $\mathcal{S}_{\mathcal{N}_j}$ be the index set of the subdomains that share the nodes in $\mathcal{N}_j$. A nec $\mathcal{N}_j$ is said to be an ancestor of another nec $\mathcal{N}_k$ if $\mathcal{S}_{\mathcal{N}_k} \subset \mathcal{S}_{\mathcal{N}_j}$. Conversely, $\mathcal{N}_k$ is an offspring of $\mathcal{N}_j$. A nec with no ancestors is called a coarse node. The RGDSW coarse basis functions are thus associated with the coarse nodes.

The value of the coarse basis functions on the interface are set as the restriction of the null space of the global Neumann problem, to each interface component $\Gamma_i$ multiplied by a partition of unity function (POU). This is required by the theory on Schwarz methods \cite{TOSELLI2005DDBOOK}. The interface operator can be written as
\begin{equation} \label{eq:rgdsw_int_bf}
    \Phi_{\Gamma}^{\text{RGDSW}} =
    \begin{pmatrix}
        R_{c_1}^{\top} \Phi_{c_1} & \dots & R_{c_n}^{\top} \Phi_{c_n}
    \end{pmatrix} =
    \begin{pmatrix}
        R_{c_1}^{\top} S_{c_1} Z_{c_1} & \dots & R_{c_n}^{\top} S_{c_n} Z_{c_n}
    \end{pmatrix},
\end{equation}
where $\Phi_{c_i}$ is the matrix whose columns are the coarse basis function associated with the coarse node $c_i$, $R_{c_i}$ is a restriction matrix from $\Gamma$ to $c_i$ and its offspring, $Z_{c_i}$ is the restriction of the null space to $c_i$ and its offspring, and $S_{c_i}$ is a scaling matrix whose diagonal entries are the evaluation of the POU. For a given node $k$, the value of the POU is given by
\begin{equation} \label{eq:rgdsw_ipou}
    \varphi (k) =
    \begin{cases}
        1 & \quad \text{if $k$ is a coarse node}, \\
        \frac{1}{|\mathcal{C}_k|} & \quad \text{otherwise},
    \end{cases}
\end{equation}
where $\mathcal{C}_k$ is the index set of coarse nodes that contain $k$ in its offspring. This is the option 1 described in \cite{DOHRMANN2017RGDSW}. Furthermore, for the Stokes flow problem studied in this work, the null space is spanned by
\begin{equation} \label{eq:null_space_stokes}
    r_{u, 1} = \begin{pmatrix}
        1,\, 0,\, 0
    \end{pmatrix}^{\top}, \;
    r_{u, 2} = \begin{pmatrix}
        0,\, 1,\, 0
    \end{pmatrix}^{\top}, \; \text{and } \,
    r_{u, 3} = \begin{pmatrix}
        0,\, 0,\, 1
    \end{pmatrix}^{\top}.
\end{equation}

% ------------

% ------------

\section{Parallel implementation} \label{sec:frosch_ngsolve}

The numerical examples discussed in the next section were implemented using the Fast and Robust Overlapping Schwarz (FROSch) package in the Trilinos toolkit \cite{HEINLEIN2020FROSCH} and the NGSolve library \cite{NGSOLVE2014}. The former implements, among others, the two-level Schwarz preconditioner described previously, while the latter provides an implementation of FE formulations, geometric modeling and mesh generation. In addition, we have developed an interface to integrate NGSolve's Python frontend to Trilinos' C++ implementation thus allowing linear systems arising from the discretization using different FE spaces to be solved with FROSch. We have used Trilinos version 16.2.0 and NGSolve version 6.2.2501. Moreover, the solution of the local subproblems and the coarse problem were computed using MUMPS 5.5.1 \cite{AMESTOY2001MUMPS}. The examples were run on the DelftBlue supercomputer \cite{DELFTBLUE2024}.

The decomposition of the interface into vertices, edges, and faces described in \cref{sec:two_level_precond} is provided by FROSch based on an algebraic procedure to identify the interface dofs. Given a non-overlapping domain decomposition, a graph that represents the mesh connectivity is obtained from the system matrix. The interface dofs are then identified by performing a breadth-first search on each subdomain to find the dofs that are not owned by the subdomain itself. Finally, the interface entities are defined based on the description in \cref{sec:rgdsw}. In our implementation, the non-overlapping domain decomposition was computed using METIS \cite{KARYPIS1998METIS}. In addition, a global indexing of the dofs is provided by NGSolve.

% ------------

% ------------

\section{Numerical results} \label{sec:results}

\begin{figure}[t]
    \centering
    \begin{subfigure}[b]{0.45\textwidth}
        \includegraphics[width=\textwidth]{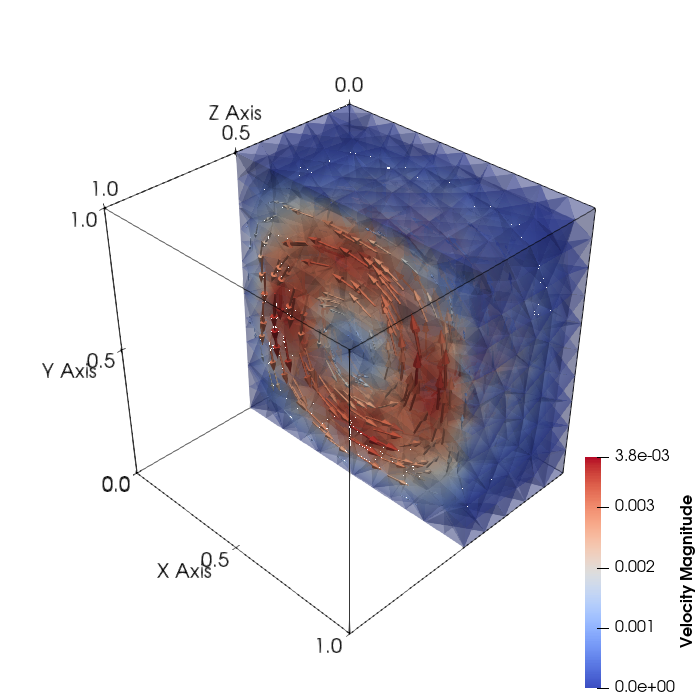}
        \caption{}
        \label{fig:solution_cube}
    \end{subfigure}
    \hfil
    \begin{subfigure}[b]{0.45\textwidth}
        \includegraphics[width=\textwidth]{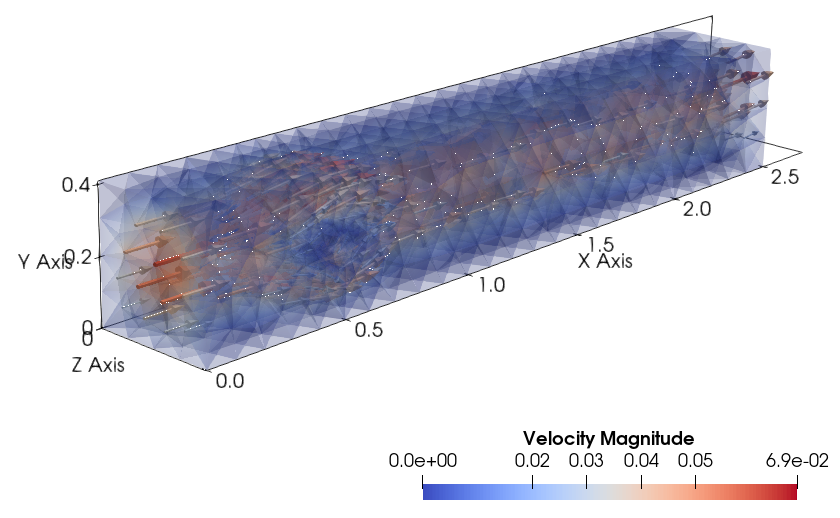}
        \caption{}
        \label{fig:solution_channel}
    \end{subfigure}
    \caption{Velocity field solutions for the proposed examples. \textbf{Left}: flow in a unit cube (slice at $x = 0.5$); \textbf{right}: flow around a cylinder in a channel.}
    \label{fig:solution_plots}
\end{figure}

\begin{figure}[t]
    \centering
    \includegraphics[width=\textwidth]{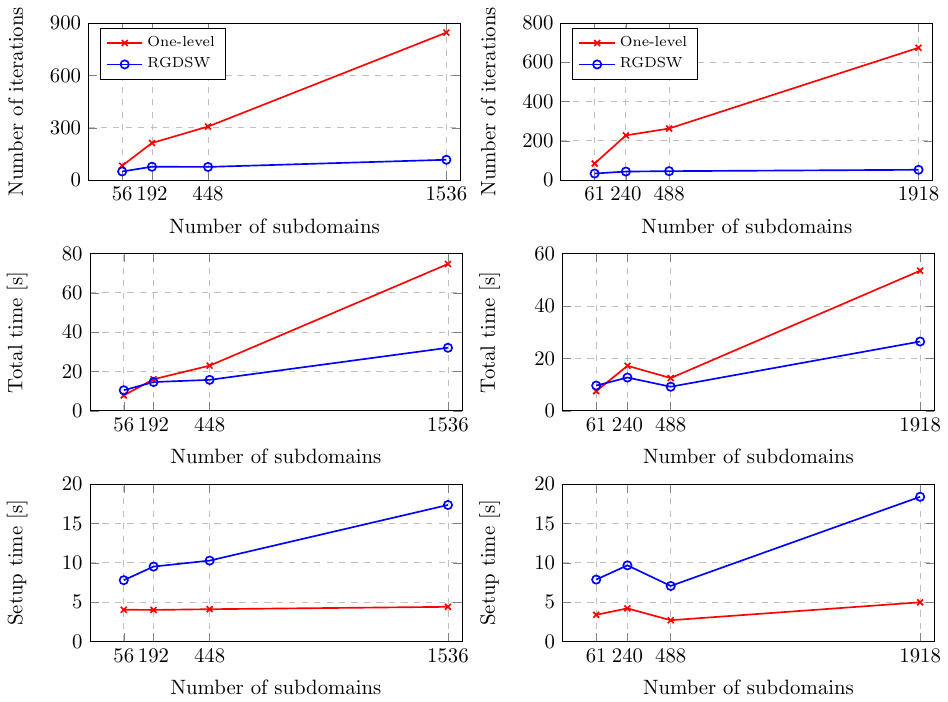}
    \caption{Numerical scalability results for the flow in a unit cube (left column) and around a cylinder (right column). \textbf{Top row:} number of iterations; \textbf{middle row:} total time to solve the system of equations (setup of the preconditioner and application of the PCG iterations); \textbf{bottom row:} time to set up the preconditioner (assembly and factorization of the local and coarse problems).}
    \label{fig:scalability}
\end{figure}

We study two numerical examples to evaluate the performance of the proposed preconditioning schemes. As a first example, we consider the Stokes equations in a unit cube $\Omega = [0, 1]^3$ with homogeneous Dirichlet boundary condition $u = (0,\, 0,\, 0)^{\top}$ and a body force term $f = (0,\, x - 1/2,\, 0)$. For the second example, we consider the flow around a cylinder with circular cross-section such that $\Omega = \Omega_{\text{channel}} \setminus \Omega_{\text{cylinder}}$, with $\Omega_{\text{channel}} = [0, 2] \times [0, 0.41]^2$ and $\Omega_{\text{cylinder}} = \{ (x,\, y,\, z) \; | \; (x - 0.5)^2 + (y - 0.2)^2 \leq 0.05^2 \wedge 0 \leq z \leq 0.41 \}$. The prescribed velocity values are
\begin{equation*}
    \begin{cases}
        \left( \frac{6 y (0.41 - y) z (0.41 - z)}{0.41^2},\, 0,\, 0 \right)^{\top} & \quad \text{on } \partial\Omega_{\text{in}} = \{ (x, y, z) \; | \; x = 0 \}, \\
        (0,\, 0,\, 0)^{\top} & \quad \text{elsewhere},
    \end{cases}
\end{equation*}
with a body force term $f = (0,\, 1,\, 0)^{\top}$. The velocity field solutions for both examples are shown in \cref{fig:solution_plots}.

The problem was discretized using the formulation in \cref{eq:bilinear_form_penalty} with a second-order Lagrangian finite element space ($P_2$) and the penalty parameter $\epsilon = 10^{-4}$. The resulting linear system was solved with the preconditioned conjugate gradient (PCG) method and the two-level OAS preconditioner described in \cref{sec:two_level_precond} with two layers of algebraic overlap. As a convergence criterion, we have adopted a reduction of the initial residual $r^{(0)} = b - A\mathbf{u}^{(0)}$ by a factor of $10^{-6}$, where $\mathbf{u}^{(0)}$ is the initial guess to the PCG iterations.

The weak scalability results for the two examples are presented in \cref{fig:scalability}. In both cases, the RGDSW preconditioner shows good numerical scalability when compared to the one-level method. The total time to solve the system of equations also scales well. The increase in the total time of the RGDSW preconditioner is mainly due to the setup of the preconditioner itself, as it can be seen in the bottom plots in \cref{fig:scalability}.

Next, we investigate the influence of the penalty term parameter $\epsilon$ on the performance of the RGDSW preconditioner. We consider the flow in a unit cube with 448 subdomains given known velocity $u_{\text{ref}}$ and pressure $p_{\text{ref}}$ solutions defined as
\begin{align*}
    \mathbf{u}_{\text{ref}} &= (x^2 (1 - x^2)^2 (2 y - 8 y^3 + 6 y^5), -y^2 (1 - y^2)^2 (2 x - 8 x^3 + 6 x^5), 0)^{\top}, \\
    p_{\text{ref}} &= x(1 - x),
\end{align*}
such that $f = -\nabla^2 \mathbf{u}_{\text{ref}} + \nabla p_{\text{ref}}$. We vary the value of $\epsilon$ and compute the $L^2$-norm of the error $||\mathbf{u}_{\text{ref}} - \mathbf{u}||_{L^2(\Omega)} = \left( \int_{\Omega} | \mathbf{u}_{\text{ref}} - \mathbf{u} | \; \text{d}x \right)^{1/2}$ in each case. The remaining solver parameters are kept the same. The results are summarized in \cref{tab:penalty_term}.

From \cref{eq:bilinear_form_penalty}, we see that, as $\epsilon \to 1$, the influence of the penalty term is reduced and the formulation approaches the type of symmetric positive definite elliptic problem for which the two-level OAS preconditioner was originally designed. Conversely, when $\epsilon \to 0$, the penalty term becomes dominant. The results in \cref{tab:penalty_term} show this effect, with an increase in the number of iterations as $\epsilon$ decreases. There is a slight increase in the total time to solve the problem, mostly due to the higher iteration count, since the setup time remains approximately constant. Moreover, the constancy of the total time allows for the parameter $\epsilon$ to be adjusted based on the desired accuracy without significantly affecting the performance of the preconditioner. An appropriate value of $\epsilon$ may thus be chosen to balance accuracy and performance.

\begin{table}[t]
\caption{Iterative performance of the RGDSW preconditioner for different values of the penalty term parameter $\epsilon$.} \label{tab:penalty_term}
\centering
\begin{tabular}{p{1.5cm}p{2.25cm}p{2.75cm}p{2cm}p{2cm}}
\hline\noalign{\smallskip}
$\epsilon$ & $||u_{\text{ref}} - u||_{L^2(\Omega)}$ & Number of iterations & Total time (s) & Setup time (s) \\
\noalign{\smallskip}\svhline\noalign{\smallskip}
$10^{-1}$ & $0.029$ & 32 & $12.27$ & $\phantom{1}9.94$ \\
$10^{-2}$ & $0.011$ & 44 & $13.19$ & $10.00$ \\
$10^{-3}$ & $0.004$ & 50 & $13.53$ & $\phantom{1}9.99$ \\
$10^{-4}$ & $0.003$ & 49 & $13.51$ & $10.02$ \\
\noalign{\smallskip}\hline\noalign{\smallskip}
\end{tabular}
\end{table}

% ------------

% ------------

% \section{Conclusions} \label{sec:conclusions}

% In this work, we have presented a scalability study of a two-level overlapping additive Schwarz preconditioner with a reduced dimension coarse space applied to the solution of the Stokes equations discretized with a reduced integration method. We have also discussed the implementation of the proposed numerical examples based on the FROSch and NGSolve libraries, allowing the iterative solution of linear systems arising from discretization using different finite element spaces with Schwarz preconditioners. Finally, our results have shown that the RGDSW coarse space has good numerical scalability for the considered examples and that the performance of the preconditioner is related to penalty parameter $\epsilon$ used in the discretization.

% ------------

\bibliographystyle{spmpsci}
\bibliography{bibliography}

@article{DOHRMANN2017RGDSW,
   author = {Clark R Dohrmann and Olof B Widlund},
   issue = {4},
   journal = {SIAM J. Sci. Comput.},
   pages = {A1466-A1488},
   title = "{On the Design of Small Coarse Spaces for Domain Decomposition Algorithms}",
   volume = {39},
   year = {2017},
}

@article{DOHRMANN2008GDSW2,
   author = {Clark R. Dohrmann and Axel Klawonn and Olof B. Widlund},
   issn = {0036-1429},
   issue = {4},
   journal = {SIAM J. Numer. Anal.},
   month = {1},
   pages = {2153-2168},
   title = "{Domain Decomposition for Less Regular Subdomains: Overlapping Schwarz in Two Dimensions}",
   volume = {46},
   year = {2008},
}

@inproceedings{DOHRMANN2008GDSW1,
   author = {Clark R. Dohrmann and Axel Klawonn and Olof Widlund},
   city = {Berlin, Heidelberg},
   isbn = {978-3-540-75199-1},
   booktitle = {Domain Decomposition Methods in Science and Engineering XVII},
   pages = {247-254},
   publisher = {Springer Berlin Heidelberg},
   title = "{A Family of Energy Minimizing Coarse Spaces for Overlapping Schwarz Preconditioners}",
   year = {2008},
}

@book{TOSELLI2005DDBOOK,
   author = {Andrea Toselli and Olof B. Widlund},
   city = {Berlin, Heidelberg},
   isbn = {978-3-540-20696-5},
   publisher = {Springer Berlin Heidelberg},
   title = "{Domain Decomposition Methods — Algorithms and Theory}",
   volume = {34},
   year = {2005},
}

@InProceedings{HEINLEIN2020FROSCH,
   author={Heinlein, Alexander and Klawonn, Axel and Rajamanickam, Sivasankaran and Rheinbach, Oliver},
   editor={Haynes, Ronald and MacLachlan, Scott and Cai, Xiao-Chuan and Halpern, Laurence and Kim, Hyea Hyun and Klawonn, Axel and Widlund, Olof},
   title={{FROSch: A Fast And Robust Overlapping Schwarz Domain Decomposition Preconditioner Based on Xpetra in Trilinos}},
   booktitle={Domain Decomposition Methods in Science and Engineering XXV},
   year={2020},
   publisher={Springer International Publishing},
   address={},
   pages={176--184},
   isbn={978-3-030-56750-7}
}

@techreport{NGSOLVE2014,
   author      = {Sch{\"o}berl, Joachim},
   title       = {{C++11 Implementation of Finite Elements in NGSolve}},
   institution = {{Institute of Analysis and Scientific Computing, TU Wien}},
   year        = {2014},
   type        = {Report},
   address     = {Vienna, Austria}
}

@misc{DELFTBLUE2024,
   author = {{D}elft {H}igh {P}erformance {C}omputing {C}entre ({DHPC})},
   title = {{D}elft{B}lue {S}upercomputer ({P}hase 2)},
   year = {2024},
   howpublished = {\url{https://www.tudelft.nl/dhpc/ark:/44463/DelftBluePhase2}},
   ark = {ark:/44463/DelftBluePhase2}
}

@article{KARYPIS1998METIS,
   author = {Karypis, George and Kumar, Vipin},
   title = {{A Fast and High Quality Multilevel Scheme for Partitioning Irregular Graphs}},
   journal = {SIAM J. Sci. Comput.},
   volume = {20},
   number = {1},
   pages = {359-392},
   year = {1998}
}

@article{AMESTOY2001MUMPS,
   author = "{Amestoy, Patrick R. and Duff, Iain S. and L'Excellent, Jean-Yves and Koster, Jacko}",
   title = "{A Fully Asynchronous Multifrontal Solver Using Distributed Dynamic Scheduling}",
   journal = "{SIAM Journal on Matrix Analysis and Applications}",
   volume = {23},
   number = {1},
   pages = {15-41},
   year = {2001},
   doi = {10.1137/S0895479899358194}
}

@InProceedings{SCHOEBERL1998,
   author="{Sch{\"o}berl, Joachim}",
   editor="{Hackbusch, Wolfgang and Wittum, Gabriel}",
   title="{Robust Multigrid Preconditioning for Parameter-Dependent Problems I: The Stokes-Type Case}",
   booktitle="{Multigrid Methods V}",
   year="1998",
   publisher="{Springer Berlin Heidelberg}",
   address="Berlin, Heidelberg",
   pages="260--275",
   isbn="978-3-642-58734-4"
}

@article{HEINLEIN2019MONOLITHICGDSW,
   author = "{Heinlein, Alexander and Hochmuth, Christian and Klawonn, Axel}",
   title = "{Monolithic Overlapping Schwarz Domain Decomposition Methods with GDSW Coarse Spaces for Incompressible Fluid Flow Problems}",
   journal = "{SIAM Journal on Scientific Computing}",
   volume = {41},
   number = {4},
   pages = {C291-C316},
   year = {2019},
   doi = {10.1137/18M1184047}
}

@article{HEINLEIN2020MONOLITHICRGDSW,
   author = "{Heinlein, Alexander and Hochmuth, Christian and Klawonn, Axel}",
   title = "{Reduced dimension GDSW coarse spaces for monolithic Schwarz domain decomposition methods for incompressible fluid flow problems}",
   journal = "{International Journal for Numerical Methods in Engineering}",
   volume = {121},
   number = {6},
   pages = {1101-1119},
   doi = {10.1002/nme.6258},
   year = {2020}
}

@article{KLAWONN2000COMPMONOOAS,
   author = "{Klawonn, Axel and Pavarino, Luca F.}",
   title = "{A comparison of overlapping Schwarz methods and block preconditioners for saddle point problems}",
   journal = "{Numerical Linear Algebra with Applications}",
   volume = {7},
   number = {1},
   pages = {1-25},
   doi = {10.1002/(SICI)1099-1506(200001/02)7:1<1::AID-NLA183>3.0.CO;2-J},
   year = {2000}
}

@article{VERFURTH1984,
   author = "{Verf\"{u}rth, R.}",
   title = "{A Multilevel Algorithm for Mixed Problems}",
   journal = "{SIAM Journal on Numerical Analysis}",
   volume = {21},
   number = {2},
   pages = {264-271},
   year = {1984},
   doi = {10.1137/0721019}
}

@article{VANKA1986,
   title = "{Block-implicit multigrid solution of Navier-Stokes equations in primitive variables}",
   author = "{S.P Vanka}",
   journal = "{Journal of Computational Physics}",
   volume = {65},
   number = {1},
   pages = {138-158},
   year = {1986},
   issn = {0021-9991},
   doi = {10.1016/0021-9991(86)90008-2}
}

@article{MALKUS1978,
    title = "{Mixed finite element methods — Reduced and selective integration techniques: A unification of concepts}",
    journal = "{Computer Methods in Applied Mechanics and Engineering}",
    volume = {15},
    number = {1},
    pages = {63-81},
    year = {1978},
    issn = {0045-7825},
    doi = {10.1016/0045-7825(78)90005-1},
    author = {David S. Malkus and Thomas J.R. Hughes}
}

\end{document}